\documentclass[12pt]{article}

\usepackage{amsmath}
\numberwithin{equation}{section}

\usepackage{amsfonts}
\usepackage{graphicx}
\usepackage[latin1]{inputenc}
\usepackage[small,margin=3mm]{caption}
\usepackage[hmargin=30mm,top=30mm,bottom=35mm]{geometry}

\newcommand{\R}{\mathbb{R}}

\newtheorem{conj}{Conjecture}[section]

\begin{document}

\title{Renormalization Group Analysis of Nonlinear Diffusion
Equations with Periodic Coefficients\footnotetext{
This work was supported by
CNPq-Brazil and
NSF under
Grant INT-0104529.
}}

\date{October, 2003}

\author{Gastão A. Braga\footnote{Departamento de Matemática - UFMG -
Brazil.} \quad
 Frederico Furtado\footnote{Department of Mathematics -
        University of Wyoming - USA.}
\\ Jussara M.
Moreira\footnotemark[1]
\quad
 Leonardo T. Rolla\footnotemark[1]
}

\maketitle

\begin{abstract}
In this paper we present an efficient numerical approach based on the
Renormalization Group method for the computation of self-similar dynamics.
The latter arise, for instance, as the long-time asymptotic behavior of
solutions to nonlinear parabolic partial differential equations.
We illustrate the approach with the verification of a conjecture about
the long-time behavior of solutions to a certain class of nonlinear
diffusion equations with periodic coefficients.
This conjecture is based on a mixed
argument involving ideas from homogenization theory and the Renormalization
Group method. Our numerical approach provides a detailed picture of the
asymptotics, including the determination of the effective or renormalized
diffusion coefficient.
\end{abstract}

This preprint has the same numbering of sections, equations and figures as the published article ``\emph{Multiscale Model. Simul. 1  (2003), p. 630--644}.''

\textbf{Key words.} renormalization group, partial differential equations, multiple scale problems,
asymptotic behavior.

\textbf{AMS subject classifications.} 35K55, 35B40, 35B33, 35B27, 74S20

\textbf{DOI.} 10.1137/S1540345902416600



\section{Introduction}
\label{sec:intro}

The time evolution of many non-equilibrium physical systems is
well described by nonlinear partial differential equations (PDEs).
It is frequently observed that under suitable, but fairly general,
conditions the time evolution of such systems becomes
asymptotically self-similar. Within the PDE framework, this
assertion amounts to the observation that the PDEs in question
have solutions that behave asymptotically as
\[
u(x,t) \sim \frac{1}{t^{\alpha}}\phi\left( \frac{x}{t^{\beta}}\right),
\,\,\,\,\,\, \mbox{as}\,\,\,\,\,\, t\rightarrow \infty.
\]
It turns out that often such a self-similar behavior possesses a great deal of
universality: the scaling exponents $\alpha$ and $\beta$ and
most aspects of the profile function $\phi(\cdot)$ are independent of initial
conditions (within a suitable class) or, more strikingly, even of
the form of the equations.

Universal behavior is a central issue in the study of
critical phenomena in equilibrium statistical mechanics and
quantum field theory. Using the Renormalization Group (RG)
approach \cite{bib:gell-low,bib:wilson1,bib:ma}, physicists
predict critical exponents and determine the universality class of
a variety of models.
In the early 90's, Goldenfeld, Oono and collaborators (see
\cite{bib:gold92} and references therein) developed a
perturbative renormalization group method for PDEs and used it
to compute perturbatively the similarity exponents in the difficult
and most frequently encountered case of self-similarity of the second kind
(in Barenblatt's classification \cite{bib:barenblatt1}).
Later, Bricmont, Kupiainen and Lin
\cite{bib:bric-kupa,bib:bric-kupa-lin} introduced a
non-perturbative renormalization group approach.
A numerical renormalization group algorithm
was developed at the same time by Chen and Goldenfeld
\cite{bib:chen}.

In this paper we combine the numerical approach of Chen and
Goldenfeld  with the constructive approach of Bricmont et al. to
develop a computationally efficient RG algorithm for the calculation
of self-similar dynamics.  This algorithm is
then used to produce a detailed picture of the asymptotic dynamics
of certain nonlinear diffusion equations with periodic coefficients;
in particular, this algorithm allows the
investigation of which aspects of these
dynamics are universal and which are not. A similar point of view
was adopted by some of the authors \cite{bib:braga-furtado-isaia} to study a
modification of Barenblatt's equation. Isaia \cite{bib:isaia02} also used
the same approach to study other classes of problems  including
some with logarithmic corrections to the time decay.

Our purpose in writing this paper is threefold.
Firstly, we conjecture the long-time asymptotic behavior
of a class of PDEs. Secondly, following the strategy of
\cite{bib:braga-furtado-isaia,bib:isaia02}, we present a systematic
numerical approach, based on the RG ideas, to study
self-similar dynamics. Finally, we point out a connection between
the renormalization group and homogenization theory
\cite{bib:papa78}.
Homogenization theory has been applied with success
to critical lattice field theories~\cite{bib:spencer97}, where
a RG approach should work. Here, in addition to adopting a
reverse point of view and showing that the RG approach
succeeds in problems typically analyzed with homogenization theory,
we further point out that the RG approach gives more general results;
see below.

Numerical procedures based on rescaling,
and thus similar in spirit to the RG approach presented here,
have previously been developed \cite{bib:papa86,bib:berger88} and
used to study solutions which blow up in finite time (as, for example,
in the case of certain reaction-diffusion equations and the cubic nonlinear
Schr\"{o}dinger equation).  Such procedures exploit the known self-similar
structure of the solutions under study to define the appropriate rescalings.
Recently, versatile and more efficient versions of such procedures
were developed \cite{bib:xiaoping2000,bib:xiaoping2003} and employed for
the computation of solutions which blow up at multiple points.
Yet, the RG procedure presented here is unique in exploiting fixed points.
The central feature of this procedure is that the rescalings are determined
dynamically, not by any a priori knowledge of the scale invariance of the
solutions.  As a result, the procedure can be used to investigate
self-similar asymptotics of both the first and second (anomalous) kind.
On the other hand the current implementation of the procedure that we
present here is not appropriate for studying blow-up problems.

This paper proceeds as follows: in Section~\ref{sec:asym}
we define the class of initial value problems we study and state
the conjecture we make about the long time
behavior of their solutions. In Section
\ref{sec:reno-grou-tran} we present a non-rigorous argument to motivate
this conjecture. In Section
\ref{sec:nume-shem} we explain how we have implemented our modification of the
numerical renormalization group scheme of Chen and Goldenfeld
\cite{bib:chen} to the oscillating equation~(\ref{eq:osci}). In
Section~\ref{sec:resu} we  present numerical results that
corroborate our conjecture. We also present studies to
validate our numerical procedure, including the computation of the critical
exponent $\alpha$ for solutions of Barenblatt's equation to illustrate the
versatility of our procedure in computing anomalous behavior. Finally, in Section
\ref{sec:conc-rema}  we present numerical evidence supporting
a similar conjecture for ``relevant" perturbations (see the Remark in the
next section), make some comments and discuss further developments.

\section{Asymptotics of oscillating equations}
\label{sec:asym}

Consider the following initial value problem (IVP):
\begin{equation}
\label{eq:osci}
u_t = [1 + \mu g(x)] u_{xx} + \lambda F(u,u_x,u_{xx})
\quad x\in\R, \; t>1,
\end{equation}
\[
u(x,1) = f(x), \,\,\,\,\,\, x\in\R.
\]
Here $\mu$ and $\lambda$ are real parameters; $\mu$ is
such that $(1 + \mu g(x))>0$ for all $x\in\mathbb{R}$; $g(x)$
is a smooth periodic function with zero mean;
$f(x)$ is a smooth and rapidly
decaying initial condition; $F(u,v,w)$ is an analytic
function of $u,v$ and $w$ around $u=v=w=0$:
\[
F(u,v,w) = u^{a}v^{b}w^{c}\quad + \sum_{\{i> a,\; j> b,\; k> c\}}
c_{ijk}u^{i}v^{j}w^{k}.
\]

Let $T$ be the period of the function $g$
and let $H(g, \mu)$ denote the harmonic
mean of $1 + \mu g(x)$:
\begin{equation}
\label{eq:harm}
H(g, \mu)\equiv
\left[ \frac{1}{T}
\int_0^{T}\frac{1}{1+\mu g(x)}\,\,\, dx
\right]^{-1}.
\end{equation}
We conjecture the following.

\begin{conj}
\label{conj:conj}
Let $u(x,t)$ be the solution of the IVP~(\ref{eq:osci}).
If the lowest order exponents of the power series expansion of
$F$ satisfy the inequality $a +2b +3c - 3 > 0$ and if the
initial condition $f$ is small enough, then
there is a constant $A$, which usually depends on $f$, $g$,
$\mu$, $\lambda$ and $F$, such that
\begin{equation}
\label{eq:long-time-beha}
t^{\alpha}u(\sqrt{t} x, t)\rightarrow
\frac{A}{\sqrt{4\pi \sigma}} e^{-\frac{x^2}{4\sigma}} \,\,\,\,\,\,
\mbox{as} \,\,\,\,\,\, t\rightarrow\infty,
\end{equation}
with $\alpha =
1/2$ and $\sigma = H(g,\mu)$.
\end{conj}

As we shall see,
this conjecture arises from formal renormalization group
arguments applied to~(\ref{eq:osci}). These arguments
are based on the methods
developed by Goldenfeld, Oono and collaborators
\cite{bib:gold92} and by Bricmont, Kupiainen and Lin
\cite{bib:bric-kupa,bib:bric-kupa-lin}. Implementing a modification
of the RG-based numerical algorithm of Chen and
Goldenfeld~\cite{bib:chen}, we verify the validity of the
conjecture.

\noindent
{\bf Remark:} Let $d_F\equiv a +2b + 3c-3$. Bricmont et al.
\cite{bib:bric-kupa-lin}
classify the perturbations $F$ (see Eq.~(\ref{eq:osci}))
as {\em irrelevant} if $d_F>0$,
{\em marginal} if $d_F=0$ ($F=uu_x, \, u_{xx}$ or $u^3$)
and {\em relevant} if $d_F<0$ ($F=u_x$ or $u^a$).
Conjecture~\ref{conj:conj} is stated for irrelevant
perturbations but similar conjectures hold for the
 marginal and relevant cases, although care is needed in order to
avoid blow up at finite time.
In particular, if $\lambda F=-u^a +$ irrelevant perturbations, $1<a<3$,
then, instead of~(\ref{eq:long-time-beha}) we expect that
\begin{equation}
\label{eq:long-time-beha-rele}
  t^{\frac 1 {a-1}}u(\sqrt{t} x, t)\rightarrow
  f_a(\frac{x}{\sqrt{\sigma}}) \,\,\,\,\,\,
  \mbox{as} \,\,\,\,\,\, t\rightarrow\infty,
\end{equation}
where $f_a$ is a function such that
$$
u(x,t) \equiv \frac{1}{t^{\frac 1 {a-1}}}f_a\left( \frac{x}{t^{1/2}}\right)
$$
is a scale invariant solution to Eq.~(\ref{eq:osci}) with $\mu=0$ and
$\lambda F=-u^a$, see~\cite{bib:brezis86,bib:bric-kupa,bib:bric-kupa-lin}.
In Section
\ref{sec:conc-rema} we present numerical results supporting
this conjecture.
Due to logarithmic corrections, our numerical approach
can not handle marginal perturbations efficiently. We point
out that logarithmic corrections have been obtained
using another version of the numerical RG in~\cite{bib:isaia02}.

We emphasize the following aspects of~(\ref{eq:long-time-beha}):
\begin{enumerate}
\item the exponent $\alpha$ is universal: $\alpha=1/2$
independently of the initial condition $f$, the periodic
function $g$, the parameters $\mu$ and $\lambda$
and the perturbation $F(u, u_x, u_{xx})$ provided
that $a+2b+3c - 3 > 0$;
\item the effective
diffusion coefficient $\sigma$ depends on $\mu$ and $g$:
$\sigma = H(g, \mu)$;
\item in general, the prefactor $A$ depends
on the initial condition $f$, on the perturbations $F(u, u_x, u_{xx})$
and $g(x)$ and on the parameters $\mu$ and $\lambda$.
\end{enumerate}

Conjecture \ref{conj:conj} can be understood as
a statement about the universality
class of Gaussian fixed points of a renormalization group transformation.
It furnishes a classification of the nonlinear terms that can be added to the
linear oscillating equation without changing the long-time asymptotics of
the associated solutions.
The work of Duro and Zuazua
\cite{bib:zuazua2000} can also be interpreted within this framework.
They have considered the equation
\[
u_t - \mbox{div}(a(\vec{x})\nabla u) = \vec{d}\cdot\nabla(|u|^{q-1}u)
\,\,\,\,\, \mbox{in} \,\,\,\,\, {\R}^N\times(0,\infty) ,
\]
where $\vec{d}\in{\R}^N$ and $a(\vec{x})$ is a smooth, symmetric
and periodic matrix. Using homogenization theory, they have
proved, among other results, the asymptotic result
(\ref{eq:long-time-beha}) whenever $q>1 + 1/N$.
Although we do not prove it, our conjecture indicates that
Duro and Zuazua's results should hold for other non-linear
equations.

\section{Renormalization group analysis}
\label{sec:reno-grou-tran}

In this section we present the formal RG arguments which led to
the aforementioned conjecture. In the sequel, we assume
that the perturbation $F(u, u_x, u_{xx})$ in Eq.~(\ref{eq:osci})
is of the form $u^au_x^bu_{xx}^c$. Later, it will become clear that
this assumption constitutes no loss of generality since the higher powers
of $F$ are ``irrelevant'' in the RG sense.
The renormalization group approach that we employ in this paper is
simply the integration of the equation followed by a rescaling.
To explain this idea, we need some preliminary notions.

Let $u$ be a real-valued function of
$(x,t) \in {\mathbb{R}} \times \mathbb{R}_+$. For a fixed
$L > 1$ and sequences of positive scaling exponents,
$\{\alpha_n\}_{n=1}^\infty$ and $\{\beta_n\}_{n=1}^\infty$, define a
sequence $\{u_n\}_{n=0}^\infty$ of rescaled functions inductively by
$u_0 = u$ and, for $n \geq 1$,
\begin{equation}
\label{eq:rescaling}
u_{n}(x,t) = L^{\alpha_{n}}u_{n-1}\left (L^{\beta_{n}}x,Lt\right ) =
L^{\alpha_{n} + \cdots + \alpha_1}u\left (L^{\beta_{n} + \cdots +
\beta_1}x,L^nt\right ).
\end{equation}
If the original function $u$ is a global solution to Eq.~(\ref{eq:osci}),
then a direct calculation reveals that $u_n$ satisfies the renormalized
initial value problem:
\begin{equation}
\label{eq:rescaled-eqb}
  \partial_t u_n = \chi_n \,
  [1 + \mu g(\omega_nx) ] \partial_x^2 u_n
+ \lambda_n \, u_n^a(\partial_x u_n)^b (\partial_x^2 u_n)^c,
\end{equation}
\[
u_n(x,1) = f_n(x).
\]
Here
$$
\chi_n = L^{n - 2 (\beta_n+\cdots +\beta_1)}, \quad
\omega_n = L^{\beta_n+\cdots+\beta_1}, \quad
$$
and
\begin{equation}
\label{eq:chin}
\lambda_n = \lambda L^{n - (b + 2c) (\beta_n+\cdots +\beta_1) + (1-a-b-c)
(\alpha_n+\cdots +\alpha_1)},
\end{equation}
and the initial data $f_n$ is
\begin{equation}
\label{eq:rescaled-eqc}
f_n(x) =
L^{\alpha_n+\cdots +\alpha_1}u\left (L^{\beta_n+\cdots +\beta_1}x,L^n\right ).
\end{equation}

The renormalization group transformation of Bricmont et al.,
which we denote as $R_L$, is now introduced.
This transformation acts on the space of initial conditions and is
defined having the linear diffusion equation
\begin{equation}
\label{eq:difu}
u_t = u_{xx}
\end{equation}
in mind.
As explained in~\cite{bib:bric-kupa,bib:bric-kupa-lin},
the study of the long time asymptotics of solutions to
(\ref{eq:osci}) with $\mu = 0$ is equivalent to
studying the fixed points, and their
basins of attraction (i.e. universality classes),
of the transformation $R_L$. The RG transformation
is designed so that its fixed points are the
similarity solutions to~(\ref{eq:difu}).
Below we adopt the same strategy and use $R_L$ to study
the long time asymptotics of solutions to the oscillating
equation~(\ref{eq:osci}).

To define $R_L$, take $L>1$ and set
$\alpha_i = \beta_i = 1/2$ for all $i\geq 1$. With this
choice, $\chi_n = 1$, $\omega_n=L^{n/2}$, $\lambda_n = \lambda L^{n[3 - (a+2b+3c)]/2}$ and the
IVP~(\ref{eq:rescaled-eqb}) becomes
\begin{equation}
\label{eq:gene-ivp}
  \partial_t u_n =
  [1 + \mu g(L^{n/2}x) ] \partial_x^2 u_n
+ \lambda L^{n[3 - (a+2b+3c)]/2}
 \, u_n^a(\partial_x u_n)^b (\partial_x^2 u_n)^c,
\end{equation}
\[
u_n(x,1) = f_n(x).
\]
If $u_n(x,t)$ is the solution to the IVP~(\ref{eq:gene-ivp})
for $t\in[1,L]$, $R_L$ is defined as
\begin{equation}
\label{eq:Rn}
(R_{L}f_n)(x) = L^{1/2}u_n\left (L^{1/2}x,L\right )\equiv f_{n+1}(x).
\end{equation}
$R_L$ depends on $n$ but we do not write explicitly this dependence
to simplify the notation.
Notice that $R_L$ has the semigroup property:
\[
R^n_Lf(x) = R_{L^n}f(x) = L^{n/2}u\left (L^{n/2}x,L^n\right ), \quad
n= 1, 2, \dots.
\]
This fact allows us to investigate the limit~(\ref{eq:long-time-beha}),
with $t=L^{n}$, by iterating $R_L$, and this is how we proceed.
We start our analysis by reviewing the results of Bricmont et al.
\cite{bib:bric-kupa,bib:bric-kupa-lin} when $\mu = 0\not= \lambda$:
\begin{enumerate}
\item let ${\cal R}$ denote the linearized RG transformation (take
$\mu = \lambda = 0$ in equation~(\ref{eq:osci})).
One can check that any multiple of the Gaussian distribution
\[
\phi_*(x)\equiv
\frac{e^{-\frac{x^2}{4}}}{\sqrt{4\pi}}
\]
is a fixed point of ${\cal R}$. $\phi_*(x)$ is called the
Gaussian fixed point;

\item any smooth function $g$ decaying sufficiently fast and
with $\hat{g}(0)=0$, where $\hat{g}(k)$ denotes the Fourier
transform of $g(x)$, contracts to zero under the action of
${\cal R}$:
\[
\|{\cal R} g\|_{\infty}\leq \frac{C}{L}\|g\|_{\infty} < \|g\|_{\infty},
\]
if $L>C$;

\item any function $f$ can be written as its projection onto $\phi_*$
plus a remainder $g$ with $\hat{g}(0)=0$:
$f = \hat{f}(0) \phi_* + g$;

\item write $R^n_Lf(x)$ as $R^n_Lf(x) = A_n \phi_*(x) + g_n(x)$,
where $\hat{g_n}(0)=0$. Then, if the exponents $a$, $b$ and $c$
are such that $[(a+2b+3c) - 3] > 0$, $A_n\rightarrow A$ and
$g_n(x)\rightarrow 0$ as $n\rightarrow \infty$.
\end{enumerate}

In other words, Bricmont et al. proved that the long-time behavior
of the solution $u(x,t)$ to the
IVP~(\ref{eq:osci}) with $\mu = 0$ is given by the solution of
\[
  u_t = u_{xx},
\]
\[
u(x,1) = A \phi_*(x).
\]

Now consider the case $\mu \not= 0 = \lambda$. Using
homogenization theory~\cite{bib:papa78},
one can conclude that the limit~(\ref{eq:long-time-beha}) exists.
This is done as follows. After $n$ RG iterations
of equation~(\ref{eq:osci}) with $\lambda=0$ we obtain:
\begin{equation}
\label{eq:line-eqb}
  \partial_t u =
  [1 + \mu g(L^{n/2}x) ] \partial_x^2 u
\end{equation}
\[
u_n(x,1) = f_n(x).
\]
Take $\epsilon_n^{-1} = L^{n/2}$. Observing that
$\epsilon_n\rightarrow 0$ as $n\rightarrow\infty$, we see that
studying the asymptotics of the solution to IVP
(\ref{eq:line-eqb}) as $n\rightarrow\infty$ is equivalent to
studying the asymptotics as $\epsilon_n\rightarrow 0$, a
situation which can be dealt with by homogenization theory which says
that the long time behavior is governed by a multiple of
the solution to
\begin{equation}
\label{eq:line-homo-eqb}
  \partial_t u =
  \sigma u_{xx}
\end{equation}
\[
u(x,1) = \phi_{\sigma}(x),
\]
where $\sigma = H(g,\mu)$ is the harmonic mean
defined by equation~(\ref{eq:harm}) and
\[
\phi_{\sigma}(x) = \frac{e^{-\frac{x^2}{4\sigma}}}{\sqrt{4\pi\sigma}}.
\]

For the general case $\mu \not= 0 \not = \lambda$, the
RG iteration of the IVP~(\ref{eq:osci}) is represented by
the IVP~(\ref{eq:gene-ivp}). As $n\rightarrow \infty$,
the second term on the right hand side of Eq.~(\ref{eq:gene-ivp})
is expected to be driven to zero if $[(a + 2b + 3c) -3]>0$.
Writing $R^n_Lf(x)$ as $R^n_Lf(x) = A_n \phi_{\sigma}(x) + g_n(x)$,
with $\hat{g_n}(0)=0$, we also expect that $A_n\rightarrow A$ and
that $g_n(x)\rightarrow 0$ as $n\rightarrow \infty$.
Therefore, for large values of $n$, we expect the solution to the
IVP~(\ref{eq:osci}) to behave as the solution to the IVP
(\ref{eq:line-homo-eqb}), with $A=\lim A_n$ and $\sigma = H(g,\mu)$.

\section{The numerical procedure}
\label{sec:nume-shem}

The RG-based numerical procedure for the integration of the IVP~(\ref{eq:osci})
is now introduced. The renormalization procedure constructs the
sequence of renormalized IVP~(\ref{eq:rescaled-eqb}) (viz. the
sequences $\{\alpha_n\}$, $\{\beta_n\}$ and $\{f_n\}$) satisfied
by the solution of~(\ref{eq:osci}) in the following steps.

Start with $f_0 = f$, the initial condition of the
IVP~(\ref{eq:osci}). For $n=0, 1, 2, \dots$
\begin{enumerate}
\item
\label{item:firststep}
  Evolve the initial data $f_n$ forward,
  from the initial time $t=1$ to $t=L$, with $L>1$ arbitrary
  but fixed, using Eq.~(\ref{eq:rescaled-eqb}), to obtain the
  solution $u_n(x,L)$.
\item
\label{item:alphastep}
  Compute the scaling exponents $\alpha_{n+1}$, $\beta_{n+1}$ as
  \begin{equation}
    \label{eq:alpha}
    L^{\alpha_{n+1}} = \frac{u_n(0,1)}{u_n(0,L)} =
    \frac{f_n(0)}{u_n(0,L)}, \quad \beta_{n+1}=\frac{1}{2}.
  \end{equation}
(This choice will be discussed shortly.)
\item
\label{item:laststep}
  Define $f_{n+1}(x) = L^{\alpha_{n+1}}u_n(L^{\beta_{n+1}}x,L)$.
\end{enumerate}

Notice that if $f_0=f$, then~(\ref{eq:rescaling}) and the definition of
$f_n$ imply that
\begin{equation}
\label{eq:fn+1}
f_{n}(x) = L^{\alpha_{n}}u_{n-1}\left (L^{\beta_{n}}x,L\right ) =
L^{\alpha_{n} + \cdots + \alpha_1}u\left (L^{\beta_{n} + \cdots +
\beta_1}x,L^n\right ),
\end{equation}
where $u(\cdot,L^n)$ is the solution of the initial value
problem~(\ref{eq:osci}) at time $t=L^n$.

A simple consideration of~(\ref{eq:fn+1}) uncovers
the relationship between the limiting behavior of $\{f_n\}$ as
$n\rightarrow\infty$ and the long time asymptotics of $u$.
To wit, rewrite~(\ref{eq:fn+1}) as
\[
u(x,L^n) = \frac{A_n}{L^{n\alpha_n}}f_n(B_n\frac{x}{L^{n\beta_n}}),
\]
where
\[
A_n = L^{n\alpha_n -
(\alpha_n+\cdots +\alpha_1)}, \quad
B_n = L^{n\beta_n - (\beta_n+\cdots +\beta_1)}.
\]
Thus, if the limits
$A_n\rightarrow A$, $B_n\rightarrow B$,
$\alpha_n\rightarrow \alpha$,
$\beta_n\rightarrow \beta$ and
$f_n\rightarrow \phi$ as $n\rightarrow\infty$
are attained, we might expect that
\[
L^{n\alpha}u(L^{n\beta}x, L^n)=
 \frac{A_n}{L^{n(\alpha_n-\alpha)}}f_n\left(B_n\frac{x}{L^{n(\beta_n-\beta)}}\right)
\rightarrow A \phi(Bx) \quad n\rightarrow \infty.
\]
This limit, in turn, establishes the self-similarity of the long
time asymptotics of $u$.

It is the numerical algorithm consisting of steps
$\ref{item:firststep}-\ref{item:laststep}$
above supplemented by the calculation of the prefactors $\{A_n\}$
and $\{B_n\}$ that we use to verify the Conjecture~\ref{conj:conj}
about the asymptotic behavior of solutions to
(\ref{eq:osci}). We remark that this procedure yields a very
detailed picture of the asymptotics. In particular, it seems ideal to
study which aspects of the dynamics are universal and which depend
on parameters and initial conditions. A number of comments are
in order.

The rationale behind the idea we use to compute $\{\alpha_n\}$
in step $\ref{item:alphastep}$ is the reputed self-similar asymptotic
dynamics we want to compute: in the self-similar regime, the solution
$u$ at $x=0$, behaves as
\[
u(0,t)\sim \mbox{const}\, t^{-\alpha}.
\]
Thus, given the relation between $u_n$ in the time interval $t\in[1,L]$
and $u$ in $t\in[L^n, L^{n+1}]$ (c.f.~(\ref{eq:rescaling})), we
expect that, as $n\rightarrow \infty$, $\alpha_n$ approaches $\alpha$.
Whereas the determination of the exponents $\alpha_n$
can always be performed as in step $\ref{item:alphastep}$ above,
the determination of the exponents $\beta_n$ is problem dependent.
It usually involves a scaling relation between the exponents
$\alpha_n$ and $\beta_n$
so that certain (a priori chosen) parts of the differential
operator remain invariant under the rescaling of step
$\ref{item:laststep}$. To illustrate this point, consider
the initial value
problem~(\ref{eq:osci}) with
$F(u, u_x, u_{xx}) = u^au_x^bu_{xx}^c$. If $\beta_n$ is chosen to be
$1/2$, then the linear operator
\begin{equation}
\label{eq:calor}
u_t = u_{xx}
\end{equation}
remains invariant under the scaling defined in step $\ref{item:laststep}$.
In this way, the dynamics associated with~(\ref{eq:calor}) can be explored.
On the other hand, if $\beta_n$ is chosen to satisfy the
scaling relation (c.f. equation~(\ref{eq:chin}))
\begin{equation}
\label{eq:scal}
1 - (b+2c)\beta_n + (1 -a -b - c)\alpha_n = 0,
\end{equation}
then it is the nonlinear operator
\begin{equation}
\label{eq:nonl}
u_t = \lambda u^au_x^bu_{xx}^c
\end{equation}
that remains invariant. In this case, it is the dynamics of
(\ref{eq:nonl}) that can be explored.
In the RG language, the choice $\beta_n = 1/2$ focuses the attention
on the dynamics of~(\ref{eq:calor}) and, as such, is suitable for
the investigation of its irrelevant and marginal perturbations.
This is done in Section~\ref{sec:resu} for the Equation~(\ref{eq:osci}).
When the perturbation is relevant, then a new scaling
relation may be needed. For instance,
with the choice~(\ref{eq:scal}), the nonlinear operator $u_t =
u^a{u_x}^b{u_{xx}}^c$ remains invariant (under the associated scaling).
Therefore this is the appropriate scaling (which then defines an
appropriate RG map) if we want to study irrelevant and marginal
perturbations to the nonlinear operator $u_t = u^a{u_x}^b{u_{xx}}^c$.
In Section~\ref{sec:conc-rema} we study a case when
both the Laplacian and the non-linear term become
marginal under an appropriate scaling (see the
remark on Section~\ref{sec:asym} and the limiting behavior
given by
(\ref{eq:long-time-beha-rele})).

Finally, the spatial scaling in step $\ref{item:laststep}$
can be realized in two different ways: by rescaling the mesh
size $\Delta x$
without changing the discrete sites $j=0, 1, \dots$, so that
after one iteration the new mesh size is
$(\Delta x)_1 = L^{-\beta_1}\Delta x$ and the new mesh points are
located at $x = jL^{-\beta_1}\Delta x$; this is the approach adopted
in Ref.~\cite{bib:chen}; or by rescaling the discrete sites
while keeping the mesh size fixed, so that after one iteration the
new discrete sites $L^{-\beta}j$ are located at $x = L^{-\beta_1}j\Delta x$.
In this approach, the values of the solution $u$ at the fixed mesh points
$x = j\Delta x$ have to be interpolated from the data given at the (new)
discrete sites in each iteration.

In our investigation of the asymptotics of problem~(\ref{eq:osci})
we have used the first approach. The successive spatial scalings performed
in step 3 change the frequency of the periodic diffusion coefficient.
It is actually this mechanism that leads to a renormalization of
the diffusion coefficient in the long time regime. Thus, to resolve
properly the effects of the oscillating diffusion, the mesh size
has to be shrunk at the same rate as the frequency is increased.

\section{Numerical results}
\label{sec:resu}

For the discretization of the governing PDE, any appropriate
scheme can be employed. We chose a simple explicit finite
difference scheme that combines Euler's method for the time discretization
with the standard three-point formula
for the discretization of the Laplacian operator and centered
differences for the first order spatial derivatives.
Given the stability constraints, the resulting scheme is
second order accurate.

To eliminate the need of numerical boundary conditions in our computations
of purely initial value problems with compactly supported initial data,
additional grid points with associted zero-valued data were added at the
boundaries of the computational domain before each timestep was performed.
In all simulations 27 grid points were used for the discretization of the
initial data in $-5 \leq x \leq 5$ (initial computational domain).
The time step was dynamically chosen to satisfy
the stability condition $(1+\max | \mu g |)
\Delta t \leq C \Delta x^2$ (the grid spacing $\Delta x$ shrinks as
the renormalization procedure is iterated).
The constant $C$ was chosen to be $0.45$.

\begin{table}
\caption{Description of the numerical simulations, with $L = 1.021$ and
$\beta=0.5$}
\begin{center} \footnotesize
\begin{tabular}{||r|r|r|r|r|r|r|r||}
\hline
\hline $n$ & $\mu$ & $g$ & $\lambda$ & $a$ & $b$ & $c$ & $f$  \\
\hline
\hline
  1 &  .1 & $g_1$   &  0  &   &   &   &  $f_1$    \\ \hline
  2 &  .1 & $g_1$   &  0  &   &   &   &  $f_2$    \\ \hline
  3 & -.15& $g_1$   &  0  &   &   &   &  $f_1$    \\ \hline
  4 &  .1 & $g_1$   &  0  &   &   &   &  $f_3$    \\ \hline
  5 &  .1 & $g_1$   & .1  & 4 & 0 & 0 &  $f_1$    \\ \hline
  6 &  .1 & $g_1$   & .1  & 2 & 1 & 0 &  $f_1$    \\ \hline
  7 &  .1 & $g_2$   & .3  & 8 & 0 & 0 &  $f_3$    \\ \hline
  8 &  .8 & $g_1$   &  0  &   &   &   &  $f_1$    \\ \hline
  9 &  .1 & $g_3$   &  0  &   &   &   &  $f_1$    \\ \hline
 10 &  .8 & $g_3$   &  0  &   &   &   &  $f_1$    \\ \hline
 11 &   0 &         &  0  &   &   &   &  $f_1$    \\ \hline
 12 &   0 &         &  0  &   &   &   &  $f_2$    \\ \hline
 13 &   0 &         &  0  &   &   &   &  $f_3$    \\ \hline
 14 &  .1 & $g_1$   & .1  & 1 & 1 & 1 &  $f_1$    \\ \hline
 15 &  .6 & $g_3$   & .1  & 0 & 1 & 1 &  $f_2$    \\ \hline
\hline
\end{tabular}
\end{center}
\label{table:table}
\end{table}

\begin{figure}[!htb]
  \begin{center}
      \includegraphics[scale=0.3]{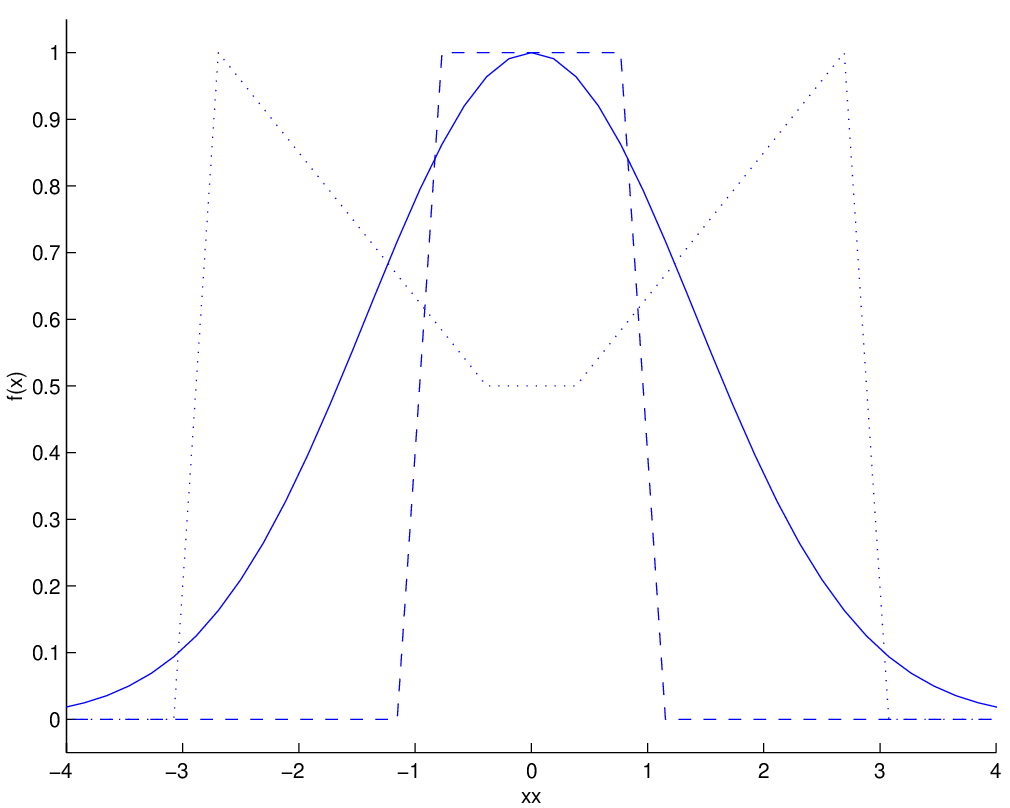}
  \end{center}
  \caption{Initial conditions $f_1$ (continuous line),
   $f_2$ (dashed line),  $f_3$ (dotted line).}
\label{fig:functions}
\end{figure}

\begin{figure}[!htb]
  \begin{center}
      \includegraphics[scale=0.3]{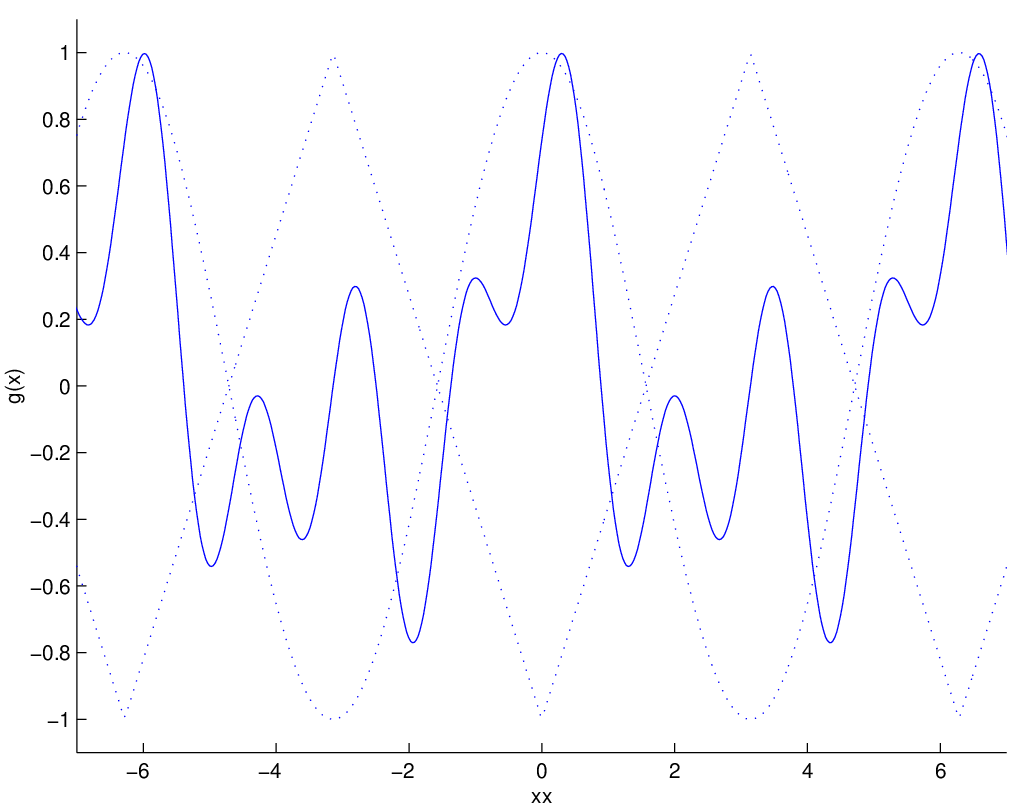}
  \end{center}
  \caption{Periodic functions: $g_1(x) = \cos(x)$, $g_2(x)=
   [ \cos(x) + \sin(2x) + \cos(4x) ]/2.72$, $g_3(x)=1-2|x/\pi {\rm mod} 2 -1|$ }
  \label{fig:perturb}
\end{figure}

Table~\ref{table:table} summarizes our numerical simulations:
the first column is the simulation number; the other columns specify
the parameters $\mu$ and $\lambda$, the exponents $a$, $b$ and $c$ and
the functions $f$ and $g$ appearing in the IVP~(\ref{eq:rescaled-eqb})
(functions $f_i$ and $g_i$, $i=1,2,3$ are specified in Figures~\ref{fig:functions} and~\ref{fig:perturb}).
In all simulations, the scaling factor
$L$ was chosen to be $L=1.021$. Although, in principle, any value of $L>1$ can be
used, for a too large value of $L$, a large number of time steps per RG iteration
would be required, and the system would shrink too quickly (due to the spatial
rescaling); for a very small $L$, a large number of RG iterations would be necessary
to produce good enough accuracy.


\begin{figure}[htb]
  \begin{center}
      \includegraphics[scale=0.3]{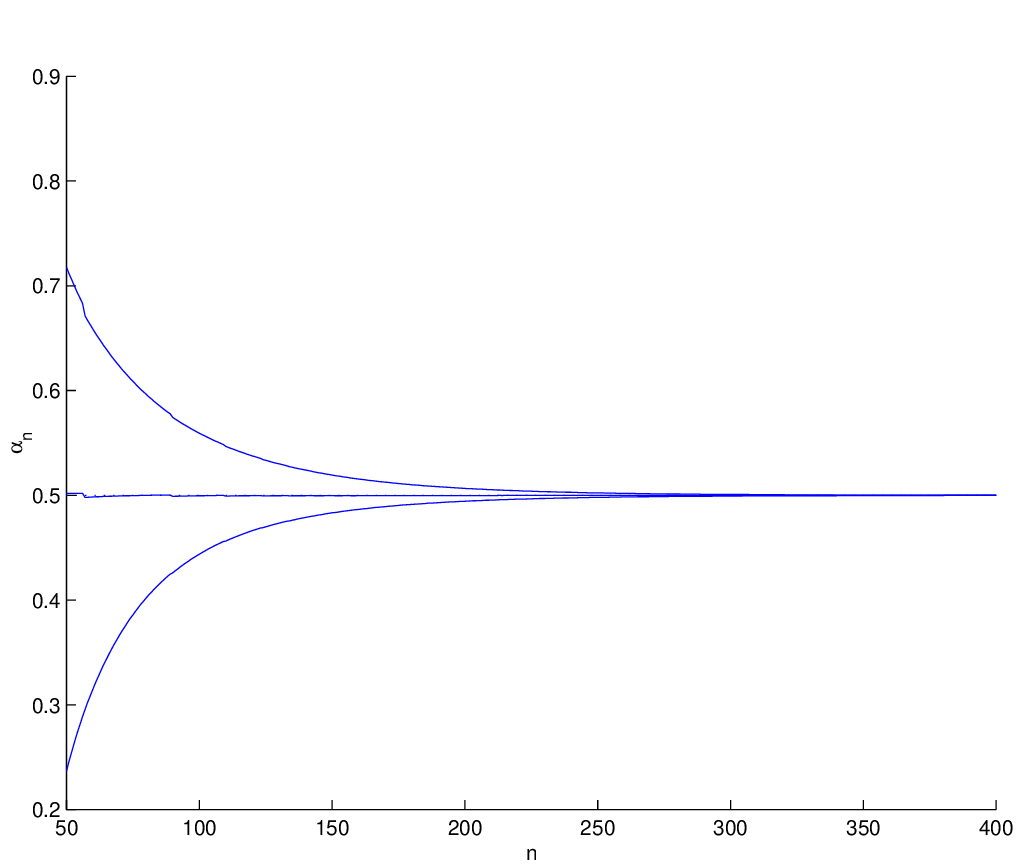}
  \end{center}
  \caption{Validation simulation (heat equation):
   plot of $\alpha_n \times n$, for
      simulations 11, 12 and 13, verifying that $\alpha^*=\lim\alpha_n$
      is universal with respect to initial conditions.}
  \label{fig:alpha-valid}
\end{figure}

\begin{figure}[!htb]
  \begin{center}
      \includegraphics[scale=0.3]{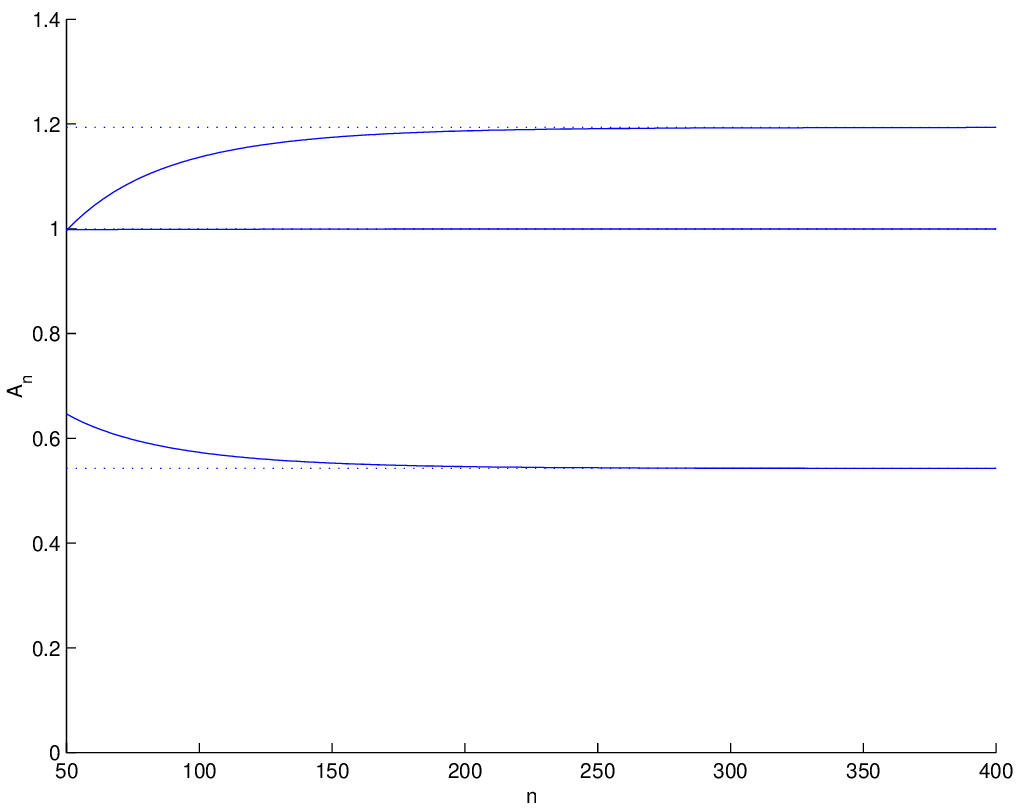}
  \end{center}
  \caption{Validation simulation (heat equation):
  plot of $A_n \times n$, for
      simulations 11, 12 and 13, verifying that $A_n\to A$ and that
      $A$ is mass dependent.}
  \label{fig:a-valid}
\end{figure}

\begin{figure}[!htb]
  \begin{center}
      \includegraphics[scale=0.3]{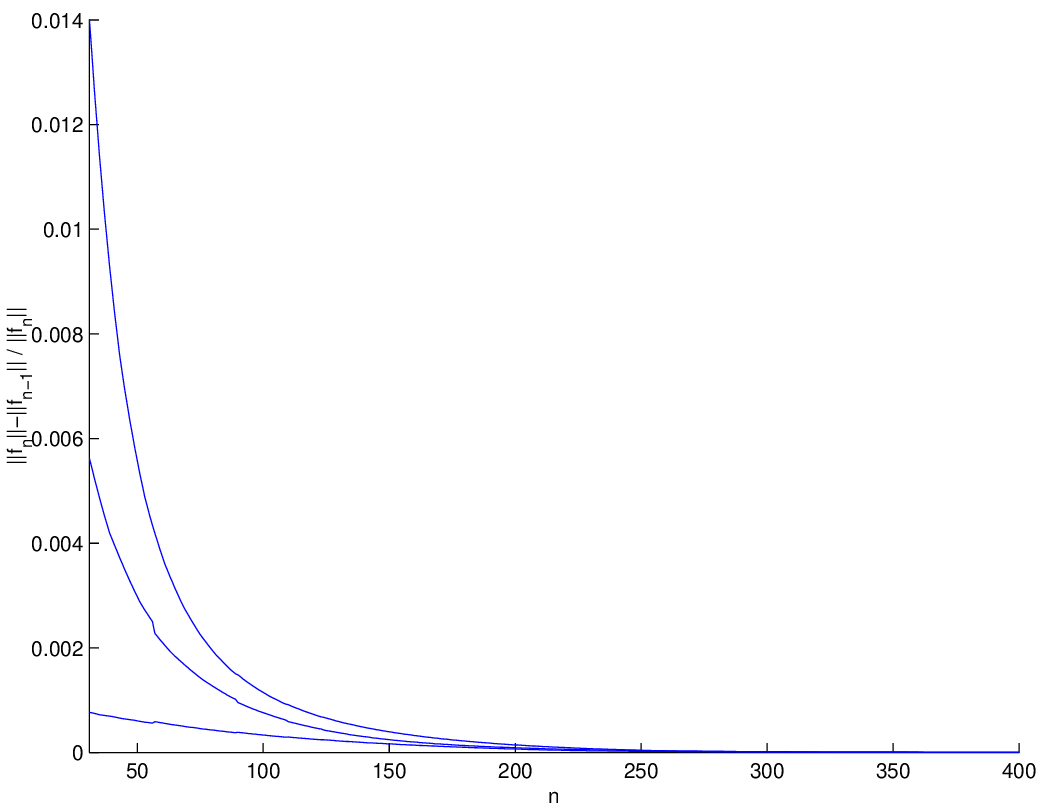}
  \end{center}
  \caption{Validation simulation (heat equation):
  plot of $\|f_n-f_{n-1}\|/\|f_n\| \times n$, for
  simulations 11, 12 and 13, verifying that the profile
    function $f_n$ converges upon the iteration of the
  numerical   RG map.}
  \label{fig:norm-valid}
\end{figure}

\begin{figure}[!htb]
  \begin{center}
      \includegraphics[scale=0.3]{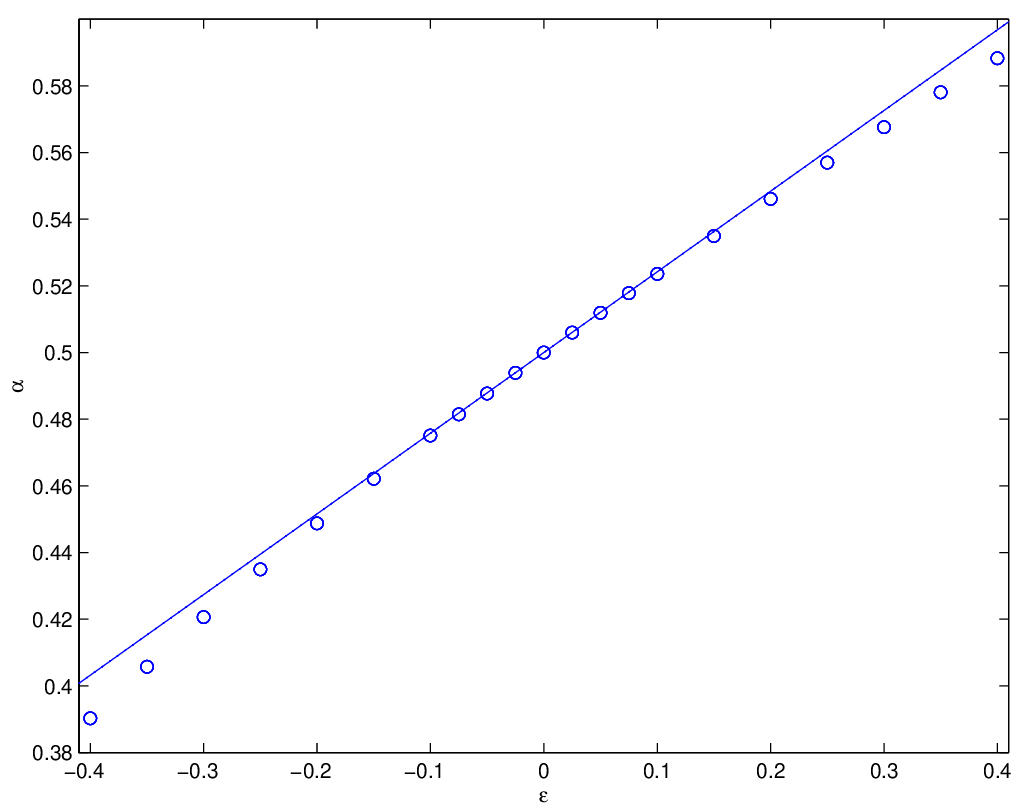}
  \end{center}
  \caption{Relation between $\alpha$ and $\epsilon$ for Barenblatt's
  equation. The continuous curve is the first-order approximation
  for $\alpha(\epsilon)$. The circle points are the results
  of our simulations.}
  \label{fig:barenblatt}
\end{figure}

We now discuss some of the simulations.
In order to validate the numerical procedure and to illustrate its
accuracy, efficiency and the amount of detailed information it provides, the
algorithm was applied to the heat equation ($\mu=0=\lambda$),
an equation whose asymptotic behavior is well known
(simulations number 11, 12, 13).
Figures~\ref{fig:alpha-valid},~\ref{fig:a-valid},~and~\ref{fig:norm-valid}
refer to this study, as explained below.
Also, simulations for  Barenblatt's equation
$$
u_t = (1 +\epsilon H(-u_t))u_{xx},
$$
where $H(u)$ is the Heaviside function, were performed
for different values of $\epsilon$.
Figure~\ref{fig:barenblatt} shows the anomalous exponent
$\alpha$ computed by our simulations,
compared with the first-order approximation~\cite{bib:gold92}
$$\alpha = \frac 1 2 + \frac \epsilon {\sqrt{2 \pi e}} + O(\epsilon^2).$$

Figure~\ref{fig:alpha-valid} depicts the convergence of the computed exponent
$\alpha_n$ to the theoretical value $\alpha = 1/2$ as the numerical RG map is iterated
for several distinct initial conditions. This illustrates the universality of
$\alpha$ with respect to sufficiently localized (e.g. with compact support)
initial conditions. Figure~\ref{fig:a-valid}
shows the convergence of the prefactors $A_n$ to the theoretical values $A$
(dotted lines). The latter depend on the mass of the initial condition.

We also verified that the profile function is equal to the Gaussian
distribution, in agreement with the analytical results.
Finally, in Figure~\ref{fig:norm-valid} we plot the relative difference
between successive profiles, $||f_n - f_{n-1}||/||f_n||$, as a function
of the number $n$ of RG iterations. The convergence of this difference
to zero can be used as a diagnosis of the self-similar behavior of the
long-time asymptotics; it can also be used to check the convergence of
the procedure and as a practical stopping criterion for the number of RG iterations.

\begin{figure}[!htb]
  \begin{center}
      \includegraphics[scale=0.3]{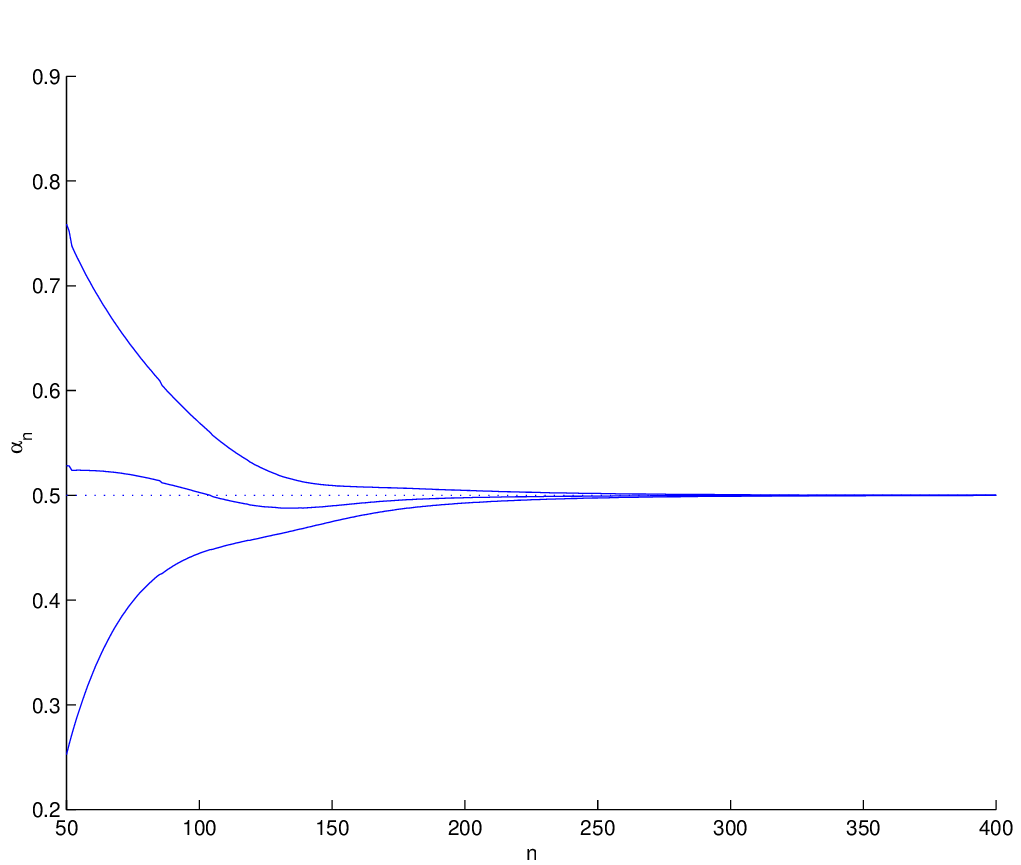}
  \end{center}
  \caption{Plot of $\alpha_n \times n$, for
      simulations 1, 2 and 4, verifying that $\alpha$ is universal with
      respect to the initial condition.}
  \label{fig:alpha-ic}
\end{figure}

\begin{figure}[!htb]
  \begin{center}
      \includegraphics[scale=0.3]{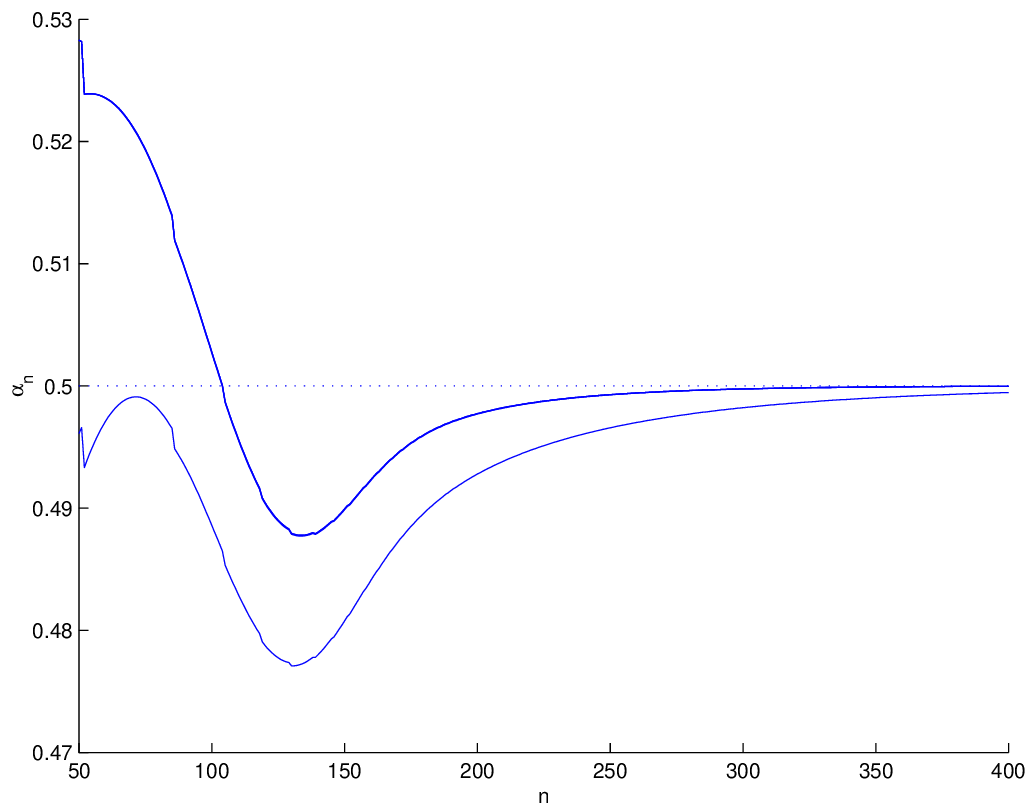}
  \end{center}
  \caption{Plot of $\alpha_n \times n$, for
      simulations 1, 5, 6, 14 and 15, verifying that $\alpha$ is universal with
      respect to the addition of nonlinear terms.}
  \label{fig:alpha-mu}
\end{figure}

\begin{figure}[!htb]
  \begin{center}
      \includegraphics[scale=0.3]{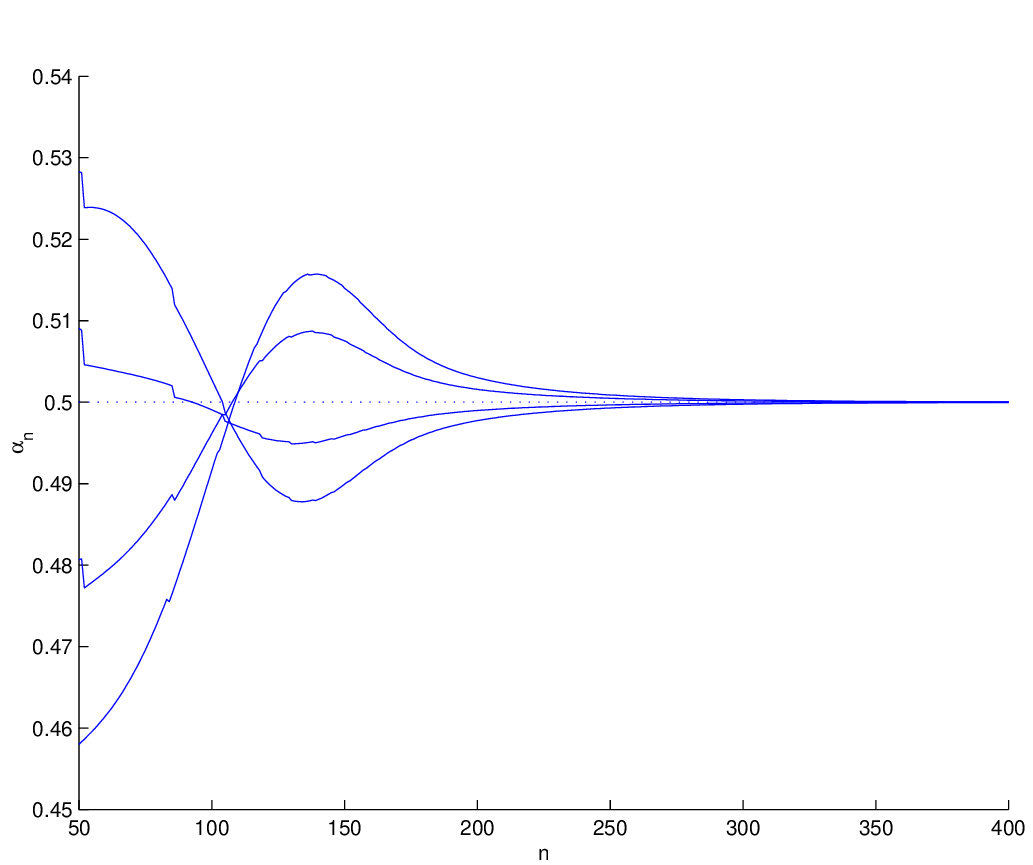}
  \end{center}
  \caption{Plot of $\alpha_n \times n$, for
      simulations 1, 3, 7 and 9, verifying that $\alpha$ is universal with
      respect to the periodic perturbation $g(\cdot)$ and the parameter
      $\mu$.}
  \label{fig:alpha-g}
\end{figure}

\begin{figure}[!htb]
  \begin{center}
      \includegraphics[scale=0.3]{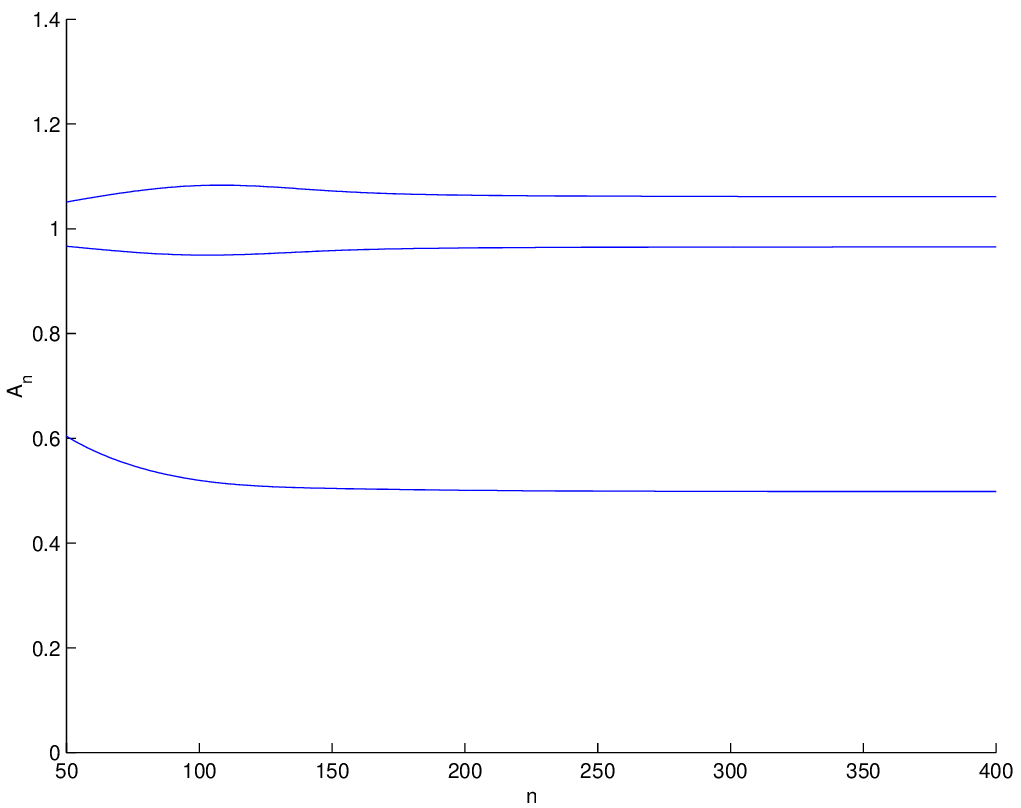}
  \end{center}
  \caption{Plot of $A_n \times n$, for simulations 1, 2 and 3,
      showing that $A_n$ does converge as $n\to\infty$.}
  \label{fig:prefactor}
\end{figure}

The remaining figures, Figures~\ref{fig:alpha-ic}-\ref{fig:prefactor},
illustrate various aspects of the asymptotics
of the solutions to the IVP~(\ref{eq:osci}) and corroborate our Conjecture
\ref{conj:conj}.
In all simulations performed, $\alpha_n$ was observed to converge quickly
to $\alpha=\frac 1 2$, as the number of RG iterations increases.
In order to verify the universality of this fact, we varied:
initial conditions (Figure~\ref{fig:alpha-ic});
nonlinear perturbations (Figure~\ref{fig:alpha-mu});
values of $\mu$ and periodic functions (Figure~\ref{fig:alpha-g}).
Figure~\ref{fig:prefactor} shows the convergence of the prefactors
$A_n$ and their nonuniversality. We emphasize that in all simulations,
the convergence of the relative difference $\|f_n-f_{n-1}\| / \|f_n\|$
(in both $L^1$ and $L^\infty$ norms)
to zero was observed, as the number of RG iterations increased, thus
providing strong evidence of the self-similar nature of the asymptotics.

\begin{figure}[!htb]
  \begin{center}
      \includegraphics[scale=0.3]{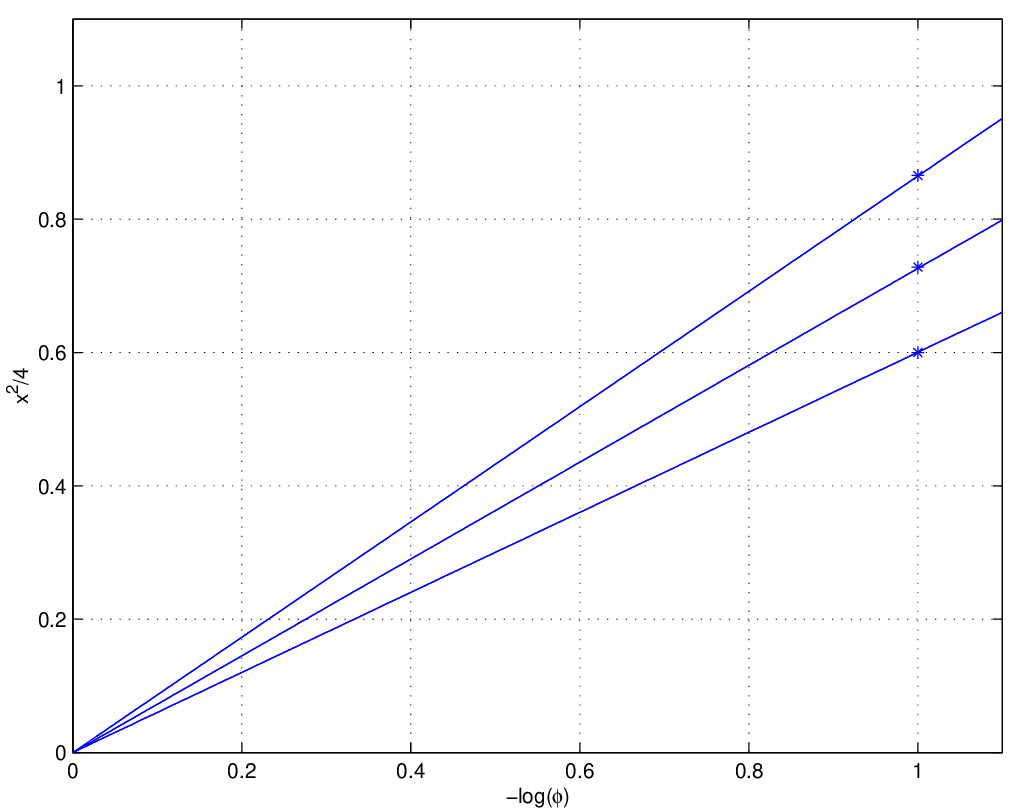}
  \end{center}
  \caption{Plot of $-\log(\phi_*) \times -\log(\phi)$, for
simulations 8, 10 and 15. The stars plotted are the
theoretically computed values for $\sigma$.}
  \label{fig:lglg}
\end{figure}

\begin{figure}[!htb]
  \begin{center}
      \includegraphics[scale=0.3]{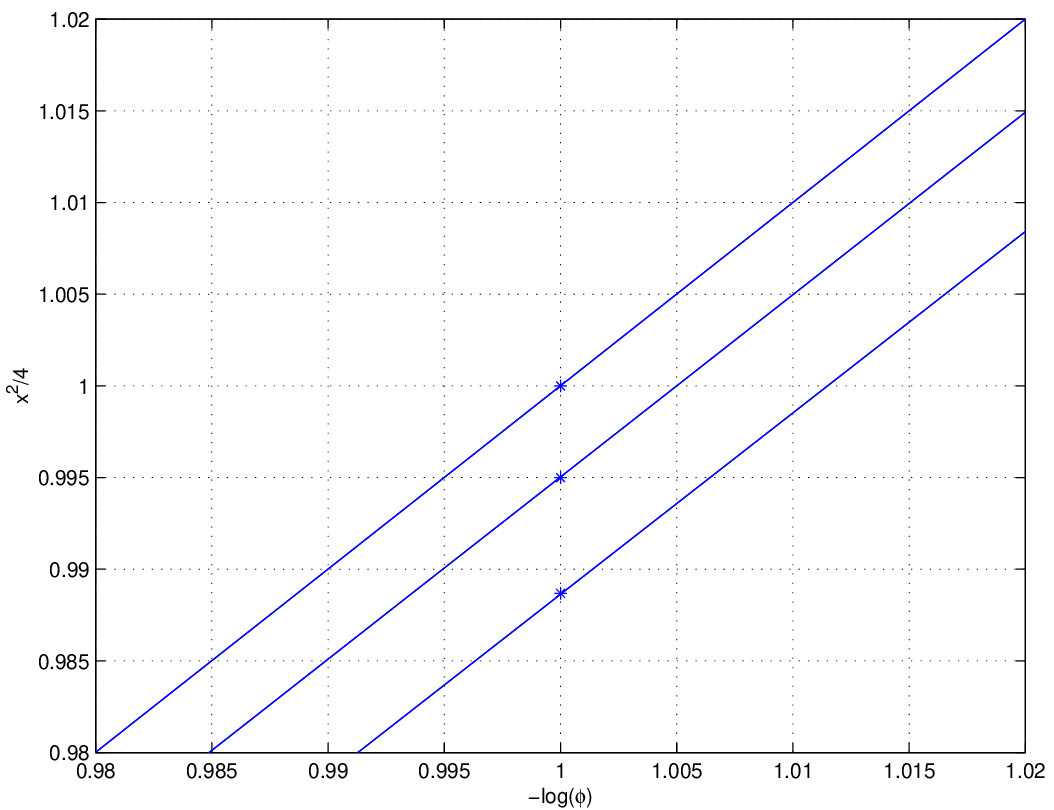}
  \end{center}
  \caption{Plot of $-\log(\phi_*) \times -\log(\phi)$, for simulations
           1, 3 and 11.
           The stars are the explicitly calculated $\sigma$.}
  \label{fig:lglgzoom}
\end{figure}

Finally, the effect of the term $\mu g$ on the diffusion coefficient
is illustrated in Figures~\ref{fig:lglg} and~\ref{fig:lglgzoom}.
These figures display the points $(-\log \phi_*, -\log \phi)$,
where $\phi$ denotes the profile function $f_n$ computed after the last RG
iteration and $\phi_*(x) = \exp(-x^2/4)$.
It is clear that for each simulation these points lie on a straight line.
This implies that the computed profiles $\phi$ are the
``renormalized Gaussians'' $\phi_\sigma(x) = \exp(-x^2/4\sigma)$,
with $\sigma$ being given by the slope of the corresponding straight line.
The theoretical values for $\sigma$ are also plotted (stars)
in these figures and they are obtained from the explicit
computation of the harmonic mean of $1+\mu g$.

\section{Concluding Remarks}
\label{sec:conc-rema}

From the numerical simulations we have a solid basis
that indicates the correctness of our Conjecture~\ref{conj:conj}.
Furthermore, using the numerical algorithm of Section
\ref{sec:nume-shem}, we have some preliminary results
for the dynamics for the following equation (which is in
divergence form):
\begin{equation}
\label{eq:modi-osci}
u_t = \partial_x ([1 + \mu g( x)] \partial_x u) +
\lambda F(u,u_x,u_{xx})
\quad x\in\R, \; t>1.
\end{equation}

Equations~(\ref{eq:osci}) and~(\ref{eq:modi-osci})
differ by a multiple of $g'( x)u_x$. Our results
indicate that this term is irrelevant in the RG sense even though
the term $u_x$ is relevant. When $F=u^qu_x$ in equation
(\ref{eq:modi-osci}), our conjecture has been
proved~\cite{bib:zuazua2000}.

The version of the numerical procedure that we described here is not appropriate
for the analysis of the dynamics when marginal perturbations of the type $uu_x$
are included. Marginal perturbations have been handled successfully using a
modified version of the procedure in~\cite{bib:isaia02}.

\begin{figure}[!htb]
  \begin{center}
      \includegraphics[scale=0.3]{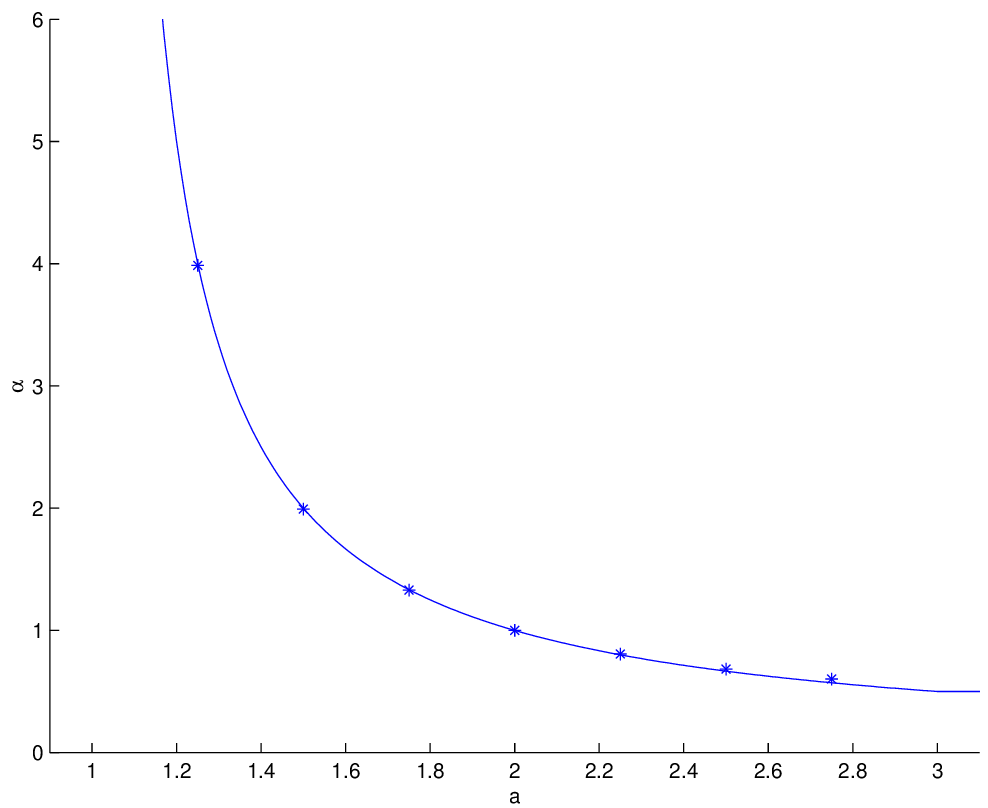}
  \end{center}
  \caption{Relation between $\alpha$ and $a$ for $\lambda F(u)=-u^a$,
  with $1<a<3$. The continuous curve is obtained theoretically
  and the stars are the results
  of our simulations.}
  \label{fig:relevant}
\end{figure}

To verify the strength and flexibility of the method,
we also performed a numerical study of
Equation~(\ref{eq:osci}) with
$\lambda F= - u^a$, $1<a<3$. First suppose
$\mu = 0$ in Eq.~(\ref{eq:osci}).
In this case, the  Laplacian and the nonlinear term
$\lambda F$ are equally important asymptotically.
The choice $\beta_n = \frac 1 2$ keeps
the Laplacian invariant to
scaling.
As discussed in
Section~\ref{sec:nume-shem},
we expect that $\alpha_n$ will converge to $\alpha = \frac 1 {a-1}$
as $n\to \infty$, so that the nonlinear
term remains in the limiting equation describing
the long-time asymptotic behavior.
This result is rigorous and it
was first proved by Brezis, Peletier and Terman~\cite{bib:brezis86}
and also
verified by Bricmont et al.~\cite{bib:bric-kupa,bib:bric-kupa-lin} using the
RG approach.
When $\mu\not = 0$, we expect that the same scalings
will hold. These scalings lead to the conjecture
(\ref{eq:long-time-beha-rele}), which we have
verified numerically.
In Figure~\ref{fig:relevant} we plot the curve
$\alpha = \frac 1 {a-1}$
together with the results of our simulations.
The renormalization of the profile function by a factor of $\sqrt{\sigma}$
in the argument has been verified:
we have plotted the computed profile function against
$x/\sqrt{\sigma}$, where the factor $\sigma$ was computed a priori
as the harmonic mean of $1+\mu g(x)$, and the plot
coincides with the plot of the computed profile function
with $\mu = 0$.


\paragraph{Acknowledgments}
L.R. and J.M. thank CNPq-Brazil
for providing them with scholarships
at the undergraduate and
graduate levels, respectively. F.F. thanks the
Departamento de Matem\'atica - UFMG for their hospitality during
his visit in July of 2002.


\begin{thebibliography}{10}

\bibitem{bib:barenblatt1}
{\sc G.~I. Barenblatt}, {\em Scaling, self-similarity and intermediate
  asymptotics}, Cambridge University Press, Cambridge, 2~ed., 1996.

\bibitem{bib:papa78}
{\sc A.~Bensousan, J.~Lions, and G.~Papanicolaou}, {\em Asymptotic Analysis of
  Periodic Structure}, North Holland, Amsterdam, 1978.

\bibitem{bib:berger88}
{\sc M.~Berger and R.~Kohn}, {\em A rescaling algorithm for the numerical
  calculation of blowing-up solutions}, Comm. Pure Appl. Math., 41 (1988),
  pp.~841--863.

\bibitem{bib:braga-furtado-isaia}
{\sc G.~A. Braga, F.~Furtado, and V.~Isaia}, {\em Smoothing of barenblatt
  equation singularity: a numerical and analytical renormalization group
  study}.
\newblock In preparation.

\bibitem{bib:brezis86}
{\sc H.~Brezis, L.~A. Peletier, and D.~Terman}, {\em A very singular solution
  of the heat equation with absorption}, Archive Rat. Mech. Anal., 95 (1986),
  pp.~185--209.

\bibitem{bib:bric-kupa}
{\sc J.~Bricmont and A.~Kupiainen}, {\em Renormalizing partial differential
  equations}, Constructive Physics, 446 (1995), pp.~83--115.

\bibitem{bib:bric-kupa-lin}
{\sc J.~Bricmont, A.~Kupiainen, and G.~Lin}, {\em Renormalization group and
  asymptotics of solutions of nonlinear parabolic equations}, Communications in
  Pure and Applied Mathematics, 47 (1994), pp.~893--922.

\bibitem{bib:chen}
{\sc L.~Chen and N.~Goldenfeld}, {\em Numerical renormalization group
  calculations for similarity solutions and travelling waves}, Physical Review
  E, 51 (1995), pp.~5577--5581.

\bibitem{bib:zuazua2000}
{\sc G.~Duro and E.~Zuazua}, {\em Large time behavior for convection-diffusion
  equations in rn with periodic coefficients}, Journal of Differential
  Equations, 167 (2000), pp.~275--315.

\bibitem{bib:xiaoping2003}
{\sc G.~Fibich, W.~Ren, and X.~P. Wang}, {\em Stability of solitary waves for
  nonlinear schrödinger equations with inhomogeneous nonlinearities}, Phys. D,
  175 (2003), pp.~96--108.

\bibitem{bib:gell-low}
{\sc M.~Gell-Mann and F.~E. Low}, {\em Quantum electrodynamics at small
  distances}, Phys. Rev., 95 (1954), pp.~1300--1312.

\bibitem{bib:gold92}
{\sc N.~Goldenfeld}, {\em Lectures on Phase Transitions and the Renormalization
  Group}, Addison-Wesley, Reading, 1992.

\bibitem{bib:isaia02}
{\sc V.~Isaia}, {\em Intermediate Asymtotic Behavior of Nonlinear Parabolic
  PDEs via a Renormalization Group Approach: a Numerical Study}, PhD thesis,
  Department of Mathematics, University of Wyoming, Laramie, Wyoming, 2002.

\bibitem{bib:ma}
{\sc S.~K. Ma}, {\em Modern Theory of Critical Phenomena}, Benjamin/Cummings
  Publishing Company, Reading, 1976.

\bibitem{bib:papa86}
{\sc D.~W. McLaughlin, G.~C. Papanicolaou, C.~Sulem, and P.~L. Sulem}, {\em
  Focusing singularity of the cubic schr{\"{o}}dinger equation}, Phys. Rev. A,
  34 (1986), pp.~1200--1210.

\bibitem{bib:spencer97}
{\sc A.~Naddaf and T.~Spencer}, {\em On homogenization and scaling limit of
  some gradient perturbations of a massless free field}, Comm. Math. Phys., 183
  (1997), pp.~55--84.

\bibitem{bib:xiaoping2000}
{\sc W.~Ren and X.~P. Wang}, {\em An iterative grid redistribution method for
  singular problems in multiple dimensions}, J. Comp. Phys., 159 (2000),
  pp.~246--273.

\bibitem{bib:wilson1}
{\sc K.~Wilson}, {\em Renormalization group and critical phenomena i, ii},
  Phys. Rev. B, 4 (1971), pp.~3174--3183, 3184--3205.

\end{thebibliography}
\end{document}